\documentclass[10pt,journal]{IEEEtran} 

\usepackage{amsmath}
\usepackage{amsfonts}
\usepackage{bm}
\usepackage{lipsum}
\usepackage{tikz}
\usepackage{textgreek}
\usepackage[normalem]{ulem}
\usepackage{soul}
\usepackage[compress]{cite}
\usepackage{wrapfig}

\renewcommand{\i}{\text{i}}

\newcommand{\stkout}[1]{\ifmmode\text{\sout{\ensuremath{#1}}}\else\sout{#1}\fi}

\title{A Local Sensory and Control Strategy for Following Hydrodynamic Signals}
\author{Brendan Colvert and Eva Kanso}
\date{}

\begin{document}
\maketitle

\begin{abstract}
Many aquatic organisms are able to track ambient flow disturbances and locate their source.  
These tasks are particularly challenging because they require the organism to sense local flow information and respond accordingly. 
Details of how these capabilities emerge from the interplay between neural control and mechano-sensory modalities remain elusive. 
Inspired by these organisms, we develop a mathematical model of a mobile sensor designed to find the source of a periodic flow disturbance. 
The sensor locally extracts the direction of propagation of the flow signal and adjusts its heading accordingly.
We show, in a simplified flow field and under certain conditions on the controller, that the mobile sensor converges unconditionally to the source of the flow field. 
Then, through carefully-conducted numerical simulations of flow past an oscillating airfoil, we assess the behavior of the mobile sensor in complex flows and demonstrate its efficacy in tracking the flow signal and locating the airfoil. 
The proposed sensory and control strategy is relevant to the design of bio-inspired underwater robots, but the general idea of orienting opposite to the direction of information propagation can be applied more broadly in optimal sensor placements and climate models.

\end{abstract}


\section{Introduction}

Fluids are essential to life.
Be it air or water, organisms at various length scales must interact with the surrounding fluid for survival.
These fluids are close to transparent, and visual sensory systems often fail to detect their motions.
Yet, fluid motions are typically characterized by long-lived coherent structures and patterns that contain information about the source generating them~\cite{Spedding2014}.
A growing body of experimental evidence suggests that organisms at various length and time scales rely on mechanosensory systems to extract such information.
For example, fish discern information contained in various flow patterns using the lateral-line sensory system~\cite{Engelmann2000}.
This ``distant-touch" capability is manifested in behaviors such as rheotaxis (alignment with or against the flow)~\cite{Arnold1974,Colvert2016}, obstacle avoidance~\cite{Windsor2010}, energy extraction from ambient vortices~\cite{Liao2003}, and identification and tracking of flow disturbances left by prey~\cite{Pohlmann2001}.
These behaviors have been reported even in the absence of visual cues, as in the extreme case of the  blind Mexican  cave fish~\cite{Montgomery2001, Windsor2010}. 
Harbor seals detect very small water currents using their long facial whiskers (or vibrissae) and could process this sensory information to track the hydrodynamic wakes of self-propelled objects~\cite{Dehnhardt2001}.
At a smaller scale, copepods have incredibly diverse sensory structures~\cite{Boxshall1997,YeNi1990} and exhibit responses to hydrodynamic signals in foraging, mating and escaping \cite{Yen1998}.

This study proposes a coupled sensory and control strategy that measures local flow information and uses these measurements to trace the flow back to its generating source.
The original inspiration is biological but the applications are much broader, ranging from robotic design to naval operations. 

The problem of identifying and following hydrodynamic cues is intimately connected to that of olfactory or chemotactic search strategies.
Chemical signals play an important role in a variety of biological behaviors, such as bacteria chemotaxis toward nutrient-rich environments~\cite{Berg1972}, homing by the Pacific salmon~\cite{Hasler2012}, foraging by seabirds~\cite{Nevitt2000}, lobsters \cite{Basil1994,Devine1982} and blue crabs \cite{Weissburg1994}, and mate-seeking by copepods and moths, which follow pheromone trails~\cite{Carde1996,Carde1997}. 

Several algorithms have been proposed for driving a mobile ``agent" towards the source of a chemical signal. 
Algorithms that follow the gradient of the signal are particularly convenient when the signal itself is sufficiently high and smooth. 
In many situations, signal detection algorithms need to be combined with gradient-based algorithms; see, for example,~\cite{Huang2017,Farrell2003}.   
The group of Naomi Leonard used gradient-based methods to cooperatively control a network of mobile sensors to climb or descend an environmental gradient; see, for example,~\cite{Bachmayer2002,Ogren2004,Biyik2008,Moore2010}. 
Another variant of the gradient-based approach was developed in the group of Miroslav Krstic, by combining a nonholonomic mobile sensor with an extremum-seeking controller that drives the sensor towards the source of a static signal field; see, for example,~\cite{CochranKrstic2009,Cochranetal2009}.  
Time-dependent signals, advected by a background flow, but in one dimension, were considered in~\cite{Hamdi2009}.
In~\cite{Vergassola2007}, an ingenious `infotaxis' approach was developed that is particularly suited to applications with time-dependent, spatially-sparse signals. 
The signal generated by the source is thought of as a message sent from the source and received by the mobile sensor through a noisy communication channel. 
The sensor reconstructs the message using a Bayesian approach conditioned on prior information gathered by the sensor.
The sensor then maneuvers in a manner that seeks to balance its need to gain more information (exploration) versus its desire to drive towards the source (exploitation).

\begin{figure*}
\centering
\includegraphics{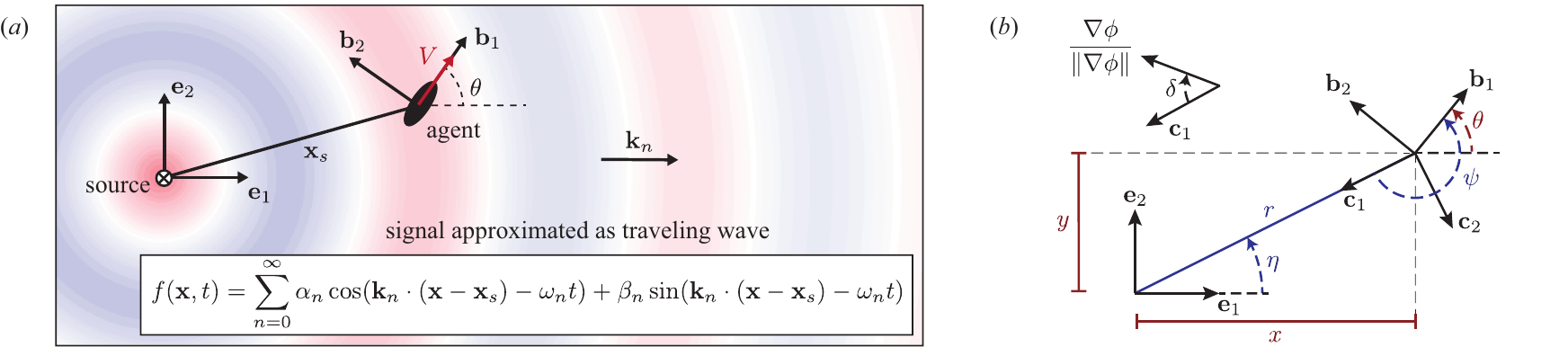}
\caption{%
(\textit{a}) A mobile sensor, located at $\mathbf{x}_s$ and equipped with body-fixed frame $(\mathbf{b}_1,\mathbf{b}_2)$, moves at 
 constant speed $V$ and changes its orientation in response to local sensory measurements of the signal field $f(\mathbf{x},t)$.
Near the sensor, $f(\mathbf{x},t)$ is expanded into a series of traveling waves with modal amplitudes $\alpha_n$ and $\beta_n$, frequencies $\omega_n$, wavenumbers $\|\mathbf{k}_n\|$ and directions $\mathbf{k}_n/\|\mathbf{k}_n\|$ indicating the local direction of signal propagation.
(\textit{b})~Relation between the inertial frame $(\mathbf{e}_1,\mathbf{e}_2)$, body-fixed frame $(\mathbf{b}_1,\mathbf{b}_2)$, where $\mathbf{b}_1$ points in the heading direction of the sensor, and body-fixed frame $(\mathbf{c}_1,\mathbf{c}_2)$, where $\mathbf{c}_1$ is in the direction pointing from the sensor to the source.
The pointing or alignment error $\delta_n$ is the angle between $\mathbf{c}_1$ and $\nabla \phi_n / \| \nabla \phi_n\|$.
}
\label{fig:setup}
\end{figure*}

In this study, we extend the notion of source seeking to problems where the relevant signal is contained in the flow field itself, as opposed to a chemical scalar field that is carried by the flow.
In a series of recent papers, we explored the information contained in various flows (e.g., shear rate, vorticity field, vortex wake patterns, etc.) and how to extract such information from local sensory measurements of the pressure and velocity fields. 
Specifically, in~\cite{Colvert2017b,Colvert2017c}, we used a physics-based approach to decipher the local flow character (extensional, rotational, or shear flow) from velocity sensors, while in~\cite{Colvert2018}, we trained a neural network to identify the global flow type (wake pattern of an oscillating airfoil) from local sensory measurements. 
Here, we develop a local sensory and control approach that is designed to operate in time-periodic signals of ``traveling wave"-like character. 
This includes the wakes generated by stationary and self-propelled bodies at intermediate Reynolds numbers.
We posit that a mobile sensor moving in the opposite direction to the direction of propagation of the signal will locate its generating source. 
We measure the signal strength locally over one period of signal oscillation and transform its time history to the frequency domain to determine the direction of propagation of the signal and its intensity. 
The mobile sensor then changes its orientation in response to this sensory information so as to orient opposite to the signal propagation. 
The coupling between the response of the mobile sensor and the signal field makes the system amenable to analysis only for relatively simple signals. 
Employing radially-symmetric signal fields and using techniques from local stability analysis and nonlinear dynamical systems theory, we rigorously prove unconditional convergence of the mobile sensor to an ``orbit-like" attractor around the source under certain assumptions on the control gain.  
To assess the behavior of the sensor in more realistic flows, we conduct careful numerical simulations of fluid flow past an oscillating airfoil and we probe the response of the mobile sensor in these more complex flow fields.


\section{Problem Statement}

We consider a mobile underwater sensor operating in an unknown hydrodynamic environment.
The mobile sensor is tasked with locating the source of a periodic flow disturbance in an unbounded two-dimensional domain.
The sensor measures the local signal contained in the fluid flow, then 
adjusts its \emph{heading} accordingly.
Our goal is to design a controller that will drive the mobile sensor to the source of the flow.

\paragraph{Mobile Sensor Model} 
We model the mobile sensor as a nonholonomic ``agent" with a sensor located at its center.
The sensor moves at a constant translational speed $V$ but is capable of changing its heading angle $\theta$ in response to its local sensory measurement (see figure~\ref{fig:setup}).
The equations of motion of the sensor are given by
\begin{equation}
\dot{\mathbf{x}}_s =  V \cos\theta \, \mathbf{e}_1 + V \sin\theta  \,\mathbf{e}_2, \qquad
\dot{\theta} = \Omega (s),
\label{eq:kinematics}
\end{equation}
where $\mathbf{x}_s$ is the position vector of the sensor with respect to an inertial frame $(\mathbf{e}_1,\mathbf{e}_2)$ and $\theta$ is measured from the positive $\mathbf{e}_1$ direction. 
The angular velocity $\Omega$ changes according to the local sensory measurement $s$.

\paragraph{Signal Field Model}
We consider an oscillating source located at the origin of the inertial frame. The source generates a hydrodynamic signal in the form of a scalar field $f(\mathbf{x},t)$ that propagates away from the source.
As the source oscillates with period $T$, we assume that the signal field reaches a steady-state limit cycle such that $f(\mathbf{x},t)=f(\mathbf{x},t-T)$ is also time periodic with period $T$.
This admits a broad class of signal fields including vortex wakes generated by biological and bioinspired engineered devices.
The exact shape of $f(\mathbf{x},t)$ is not known to the sensor.
Our goal is to steer the mobile sensor towards the origin using only time history of the signal at the sensor location.

\begin{figure*}
\centering
\includegraphics{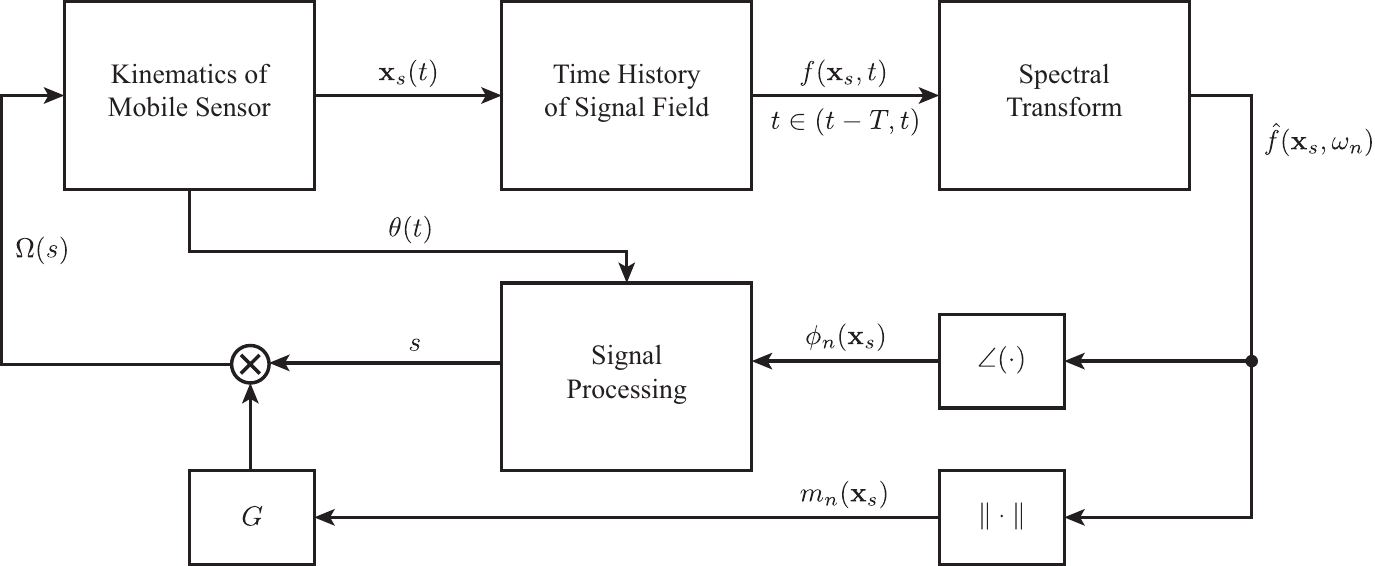}
\caption{%
Block diagram depicting the sensory and control laws of the mobile sensor.
The time history of the signal field $f(\mathbf{x}_s,t)$ is measured locally at $\mathbf{x}_s$ over one oscillation period of the signal and transformed to the frequency domain $\hat{f}(\mathbf{x}_s,\omega_n)$.
The magnitude $m_n(\mathbf{x})$ and phase fields $\phi_n(\mathbf{x})$ are computed from the spectral information to obtain the sensory output $s$ and determine the gain function $G$ of the feedback controller 
$\Omega = G s$. 
}
\label{fig:signal_diagram}
\end{figure*}

In light of the assumed structure of the signal field, we can locally (in the vicinity of the sensor at $\mathbf{x}_s$) decompose $f(\mathbf{x},t)$ as a series of traveling waves
\begin{multline}
f(\mathbf{x},t) = \sum_{n = 0}^\infty \alpha_n(\mathbf{x}) \cos \left( \mathbf{k}_n(\mathbf{x}) \cdot (\mathbf{x}-\mathbf{x}_s) - \omega_n t\right)\\
+  \beta_n (\mathbf{x}) \sin \left( \mathbf{k}_n (\mathbf{x}) \cdot (\mathbf{x}- \mathbf{x}_s) - \omega_n t \right),
\label{eq:traveling_wave}
\end{multline}
where $n$ is a nonnegative integer indicating the mode index, $\alpha_n(\mathbf{x})$ and $\beta_n(\mathbf{x})$ are the amplitudes, $\omega_n$ are the frequencies, and $\mathbf{k}_n(\mathbf{x})$ are the wavenumber vectors of the $n^{\text{th}}$ modes.
Next, we establish a framework that allows the mobile sensor to adjust its heading direction relative to the direction of propagation of the signal in order to locate the source position.
This framework relies on the premises that the signal field propagates in a traveling wave fashion and that the source-seeking agent is capable of deciphering salient information contained in these traveling waves, as depicted in figure~\ref{fig:signal_diagram} and explained next.

\paragraph{Local Sensing Model}

We define the frequency spectrum $\hat{f}(\mathbf{x},\omega)$ of the signal field $f(\mathbf{x},t)$ as follows
\begin{equation}
\hat{f}(\mathbf{x},\omega) = \lim_{\stackrel{\tau \rightarrow \infty}{\xi \rightarrow \omega}} \frac{1}{2\tau} \int_{-\tau}^\tau f(\mathbf{x},t) \exp (-\i \xi t) dt,
\label{eq:spectrum_definition}
\end{equation}
where $\i = \sqrt{-1}$ is the imaginary unit and $\omega$ is the frequency at which the spectrum is being evaluated.
We substitute~\eqref{eq:traveling_wave} into~\eqref{eq:spectrum_definition},  recalling Euler's trigonometric identities $ \cos u = \left[\exp( \i u) +  \exp (- \i u) \right]/2$ and $ \sin u = \left[\exp (\i u) -  \exp (- \i u )\right]/2\i$. 
We get that at $ \omega = \omega_n$,
\begin{equation}
\hat{f}(\mathbf{x},\omega_n) = \frac{1}{2}\sum_{n=0}^\infty (\alpha_n+\i \beta_n) \exp(-\i \mathbf{k}_n \cdot (\mathbf{x} - \mathbf{x}_s)).
\label{eq:spectrum_eval_wave}
\end{equation}
Efficient and practical estimation of the relevant frequencies $\omega_n$ contained within a signal is the subject of a great deal of established literature (see, e.g.,~\cite{Stoica2005}).
Here, we remain agnostic to the method of spectral analysis and we assume that the period $T$ of the signal $f(\mathbf{x},t)$ can be estimated offline and that $\omega_n=2\pi n/T$.

For each mode $n$, we define the spectral magnitude field $m_n(\mathbf{x})$ and the spectral phase field $\phi_n(\mathbf{x})$,
\begin{equation}
m_n(\mathbf{x}) = \| \hat{f}(\mathbf{x},\omega_n) \| \quad \text{and} \quad \phi_n(\mathbf{x}) = \angle \hat{f}(\mathbf{x},\omega_n),
\label{eq:mag_phase_definition}
\end{equation}
where $\| \cdot \|$ and $\angle (\cdot)$ denote the magnitude and argument of a complex number, respectively.
In light of~\eqref{eq:spectrum_eval_wave}, these fields are given by
\begin{equation}
m_n(\mathbf{x}) = \frac{1}{2} \sqrt{\alpha_n^2(\mathbf{x})+ \beta_n^2(\mathbf{x})},
\end{equation}
and
\begin{equation}
\phi_n(\mathbf{x}) = \tan^{-1} \frac{\beta_n(\mathbf{x})}{\alpha_n(\mathbf{x})} -\mathbf{k}_n(\mathbf{x}) \cdot (\mathbf{x}-\mathbf{x}_s) \ \text{mod} \ 2\pi .
\label{eq:phase_field_full}
\end{equation}
We expand the phase field $\phi(\mathbf{x})$ about $\mathbf{x}_s$ using Taylor series expansion,
\begin{equation}
\phi_n(\mathbf{x}) = \phi_n(\mathbf{x}_s) + \nabla \phi_n(\mathbf{x}) \bigr|_{\mathbf{x}_s} \cdot (\mathbf{x}-\mathbf{x}_s) + \mathcal{O}(\| \mathbf{x}-\mathbf{x}_s\|^2),
\label{eq:taylor_expansion}
\end{equation}
where $\nabla (\cdot)\bigr|_{\mathbf{x}_s}$ denotes the spatial gradient evaluated at $\mathbf{x}_s$.
Comparing~\eqref{eq:phase_field_full} and~\eqref{eq:taylor_expansion} and neglecting higher order terms, we obtain that 
\begin{equation}
\phi_n(\mathbf{x}_s ) \approx \tan^{-1} \frac{\beta_n}{\alpha_n} \quad \quad \text{and} \quad \quad \nabla \phi_n ({\mathbf{x}_s} )\approx -\mathbf{k}_n.
\label{eq:gradient_approx}
\end{equation}
The wavenumber vector $\mathbf{k}_n$ describes the direction of propagation of the signal in~\eqref{eq:traveling_wave}. Therefore, if the phase field $\phi_n$ and its gradient $\nabla \phi_n$ can be measured, the latter can be used as a tool to determine the  direction of signal propagation.
With this in mind, we define the \emph{spectral direction field} $\nabla \phi_n /\|\nabla \phi_n\|$, which is a unit vector field describing the local direction of signal propagation for the $n^\text{th}$ mode.

A few remarks regarding the practical computation of the phase gradient are in order. 
The assumption of direct access to the gradient is feasible because one can use multiple sensors on the sensory platform to produce a finite-difference approximation of the gradient.
A more subtle issue lies in the fact that $\phi_n(\mathbf{x})$ is a scalar field in the frequency domain whereas our mobile sensor lives in the time domain and needs to respond to the measured signal in real time. 
Measuring $\nabla \phi_n(\mathbf{x}_s)$ requires access to the time history of the $f(\mathbf{x},t)$ over one period of oscillations. 
Here, we assume that the mobile sensor moves at a much slower time scale compared to the time scale of the signal oscillations.
This translates to a constraint on the speed $V$ of the mobile sensor such that $V \ll  \omega_n / \| \mathbf{k}_n \|$ where $\omega_n / \| \mathbf{k}_n \|$ is the wavespeed of mode $n$. In this regime, we can use the quasi-steady assumption that the sensor measures a full time-period of the signal field ``instantaneously" relative to its swimming speed.

\paragraph{Sensory information and control law}
Since the first mode $n = 1$ in the ambient signal  is likely to be among the most energetic modes and therefore the most easily measured with accuracy in practical applications, we assume that it is sufficient for the sensor to measure the first mode only. Hereafter, we drop the index $n$ with the understanding that $n=1$.

We affix a body frame $(\mathbf{b}_1, \mathbf{b}_2)$ to the agent such that $\mathbf{b}_1$ is pointing in the heading direction ($\mathbf{b}_1\cdot\mathbf{e}_1= \cos \theta$ as shown in figure~\ref{fig:setup}). We consider the measurable sensory signal $s$ to be the lateral component of the phase gradient at the position $\mathbf{x}_s$ of the sensor, namely, 
\begin{equation}
s =  \dfrac{\nabla\phi({\mathbf{x}_s})}{\|\nabla\phi({\mathbf{x}_s})\|}\cdot \mathbf{b}_2 (\theta).
\end{equation}
We introduce the feedback controller $\Omega = G(m(\mathbf{x}_s)) s(\mathbf{x}_s,\theta)$,
where  $G(m(\mathbf{x}_s))$  is a proportional gain function that may dependent on the spectral magnitude field $m(\mathbf{x}_s)$ evaluated at the sensor location.
The objective of this controller is to steer the mobile sensor to follow the local gradient of the spectral phase field.
The details of both the sensory system and control model are summarized schematically in figure~\ref{fig:signal_diagram}.

\section{Equations of Motion}
We substitute the control law $\Omega = G(m(\mathbf{x}_s)) s(\mathbf{x}_s,\theta)$ into equations~\eqref{eq:kinematics} to obtain a closed system of equations that govern the motion of the mobile sensor in response to the local sensory signal $s$ it perceives,
\begin{equation}
\begin{split}
\label{eq:eom_body}
\dot{\mathbf{x}}_s = V \mathbf{b}_1(\theta),  \qquad \dot{\theta}  =  G \dfrac{\nabla\phi(\mathbf{x}_s)}{\|\nabla\phi(\mathbf{x}_s)\|} \cdot\mathbf{b}_2(\theta).
\end{split}
\end{equation}
In component form, one gets three scalar, first-order  differential equations,
\begin{equation}
\begin{split}
\label{eq:eom_scalar}
\dot{x} & = V \cos\theta, \qquad \dot{y} = V \sin\theta, \\[1ex]
\dot{\theta} & =  G \left[ \dfrac{-\dfrac{\partial \phi}{\partial x}\sin\theta + \dfrac{\partial \phi}{\partial y}\cos\theta}
{\sqrt{\left(\dfrac{\partial \phi}{\partial x} \right)^2+ \left(\dfrac{\partial \phi}{\partial y}\right)^2}} \right].
\end{split}
\end{equation}
Here, the spectral magnitude $m$ and phase $\phi$ of the signal (as well as $\partial\phi/\partial x$ and $\partial\phi/\partial y$) are all evaluated at the location $\mathbf{x}_s \equiv (x,y)$ of the sensor. For a given choice of $G(m)$, these coupled equations, together with the local sensory measurements of the signal,  constitute a closed-form system that can be solved to obtain the trajectory of the mobile sensor $(x(t),y(t))$ and its heading angle $\theta(t)$.

To analyze the ability of the mobile sensor  to track the signal field and locate its source, it is convenient to introduce a right-handed orthonormal frame $(\mathbf{c}_1,\mathbf{c}_2)$ that is attached to the sensor and pointing towards the source such that $\mathbf{c}_1 =  -\mathbf{x}/\| \mathbf{x}\|$.  It is also convenient to define the \emph{alignment error} $\delta(\mathbf{x})$ as the angle measured from $\mathbf{c}_1$ (the unit vector pointing from the sensor directly to the source) to the  direction of propagation of the signal $\nabla \phi(\mathbf{x})/\|\nabla \phi(\mathbf{x)}\|$ such that
\begin{equation}
\dfrac{\nabla \phi}{\| \nabla \phi \|}  = \cos \delta \, \mathbf{c}_1 +\sin \delta \, \mathbf{c}_2.
\label{eq:pointing_error_def}
\end{equation}
The alignment error field $\delta(\mathbf{x})$ is a local measure of the degree of alignment of the spectral direction field ${\nabla \phi}/{\| \nabla \phi \|}$ with the direction towards the source.
When $\delta = 0$,  $\nabla \phi$ points directly toward the source. For $\delta < 0$ or $\delta > 0$, $\nabla \phi$ points to the left or to the right of the source, respectively.
Since we are driving the mobile sensor to follow the local spectral direction field under the premise that it will lead to the source, the alignment error provides a measure of how well the mobile sensor will head towards the source from a given location $\mathbf{x}$ in the signal field.

We now rewrite the equations of motion~\eqref{eq:eom_body} in the $(\mathbf{c}_1,\mathbf{c}_2)$ frame. To this end, note that $(\mathbf{c}_1,\mathbf{c}_2)$ is related to 
 $(\mathbf{b}_1,\mathbf{b}_2)$  via the orthogonal transformation 
\begin{equation}
\label{eq:cb}
\left(\begin{array}{c} \mathbf{c}_1 \\ \mathbf{c}_2 \end{array} \right) =
 \left(\begin{array}{cc} \cos\psi &  \sin\psi \\ -\sin\psi  & \cos\psi \end{array} \right)
\left(\begin{array}{c} \mathbf{b}_1 \\ \mathbf{b}_2 \end{array} \right) ,
\end{equation}
and  to the inertial frame $(\mathbf{e}_1,\mathbf{e}_2)$  via the transformation
\begin{equation}
\label{eq:ce}
\left(\begin{array}{c} \mathbf{c}_1 \\ \mathbf{c}_2 \end{array} \right) =
 \left(\begin{array}{cc} - \cos\eta &  - \sin\eta \\ \sin\eta  & -\cos\eta \end{array} \right)
\left(\begin{array}{c} \mathbf{e}_1 \\ \mathbf{e}_2 \end{array} \right) .
\end{equation}
Here, the angle $\eta$ is given by the polar transformation from $(x,y)$ to $(r,\eta)$
\begin{equation}
x = r \cos\eta, \qquad y = r \sin\eta,
\end{equation}
with $r^2 = x^2 + y^2$. 
The angle $\psi$ is given by $\psi =\pi - (\theta- \eta) $. In figure~\ref{fig:setup}(b) is a depiction of the all three reference frames  $(\mathbf{e}_1,\mathbf{e}_2)$,  $(\mathbf{b}_1,\mathbf{b}_2)$ and $(\mathbf{c}_1,\mathbf{c}_2)$, together with the polar coordinate angle $\eta$ and the angles $\theta$ and $\psi$.  It is important to emphasize here that $\eta$ does not contain any information about the heading direction of the sensor but $\psi$ does.

By definition of  $(r,\eta)$ and $(\mathbf{c}_1,\mathbf{c}_2)$, we get that $\mathbf{x}_s = - r\mathbf{c}_1$ and $\dot{\mathbf{x}}_s = - \dot{r}\mathbf{c}_1 - r\dot{\eta}\mathbf{c}_2$. Substituting these expressions into~\eqref{eq:eom_body} and using~\eqref{eq:pointing_error_def} and~\eqref{eq:cb}, we get that
\begin{equation}
\begin{split}
\label{eq:eom_cframe}
\dot{r}  & = - V \cos \psi, \qquad r\dot{\eta}   = V \sin\psi \\[1ex]
r \dot{\psi} & = V \sin\psi  - rG\left( \cos \delta \sin \psi + \sin \delta \cos \psi\right).
 \end{split}
\end{equation}
Here, the alignment error $\delta(r,\eta)$ is expressed as a function of the polar coordinates $(r,\eta)$. This is rooted in the fact that the signal field $f$, and consequently its magnitude $m$ and phase $\phi$, are functions of space only, which can be expressed as functions of $(r,\eta)$ only.

By choice, we will consider three types of gain functions: (i) static gain of constant value $G_o$, (ii)  gain proportional to the  magnitude of the signal $G = G_o m(r,\eta)$, and (iii)  inversely-proportional gain  $G = G_o/m(r,\eta)$. It will be convenient for later to introduce the constant dimensionless parameter $\rho = V/G_o$, which can be thought of as a \emph{gain length scale}. In all three cases, the gain $G$ is function of the sensor position only.

The coupled equations in~\eqref{eq:eom_cframe} together with the local sensory measurement of the signal field and the choice of the gain function  constitute a closed-form system that describe the motion of the sensor in terms of $(r,\eta, \psi)$. They are equivalent to equations~\eqref{eq:eom_scalar}, except at the origin.

\section{Sensor in Simplified Signal Field} 
\label{sec:simplefield}

\begin{figure*}
\centering
\includegraphics{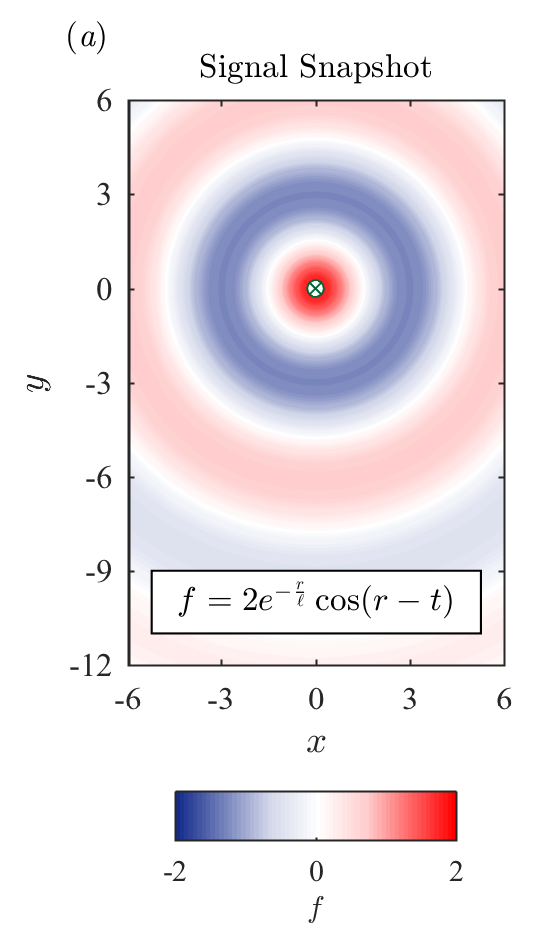}
\includegraphics{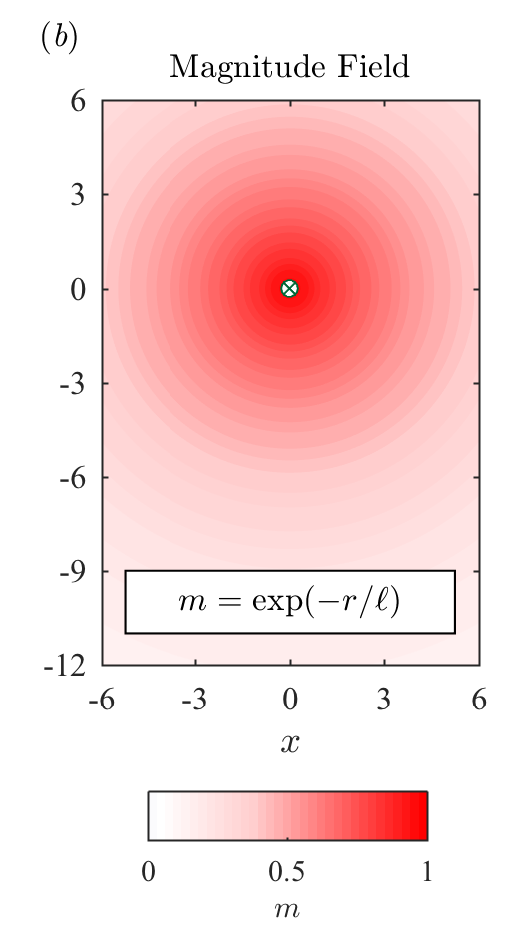}
\includegraphics{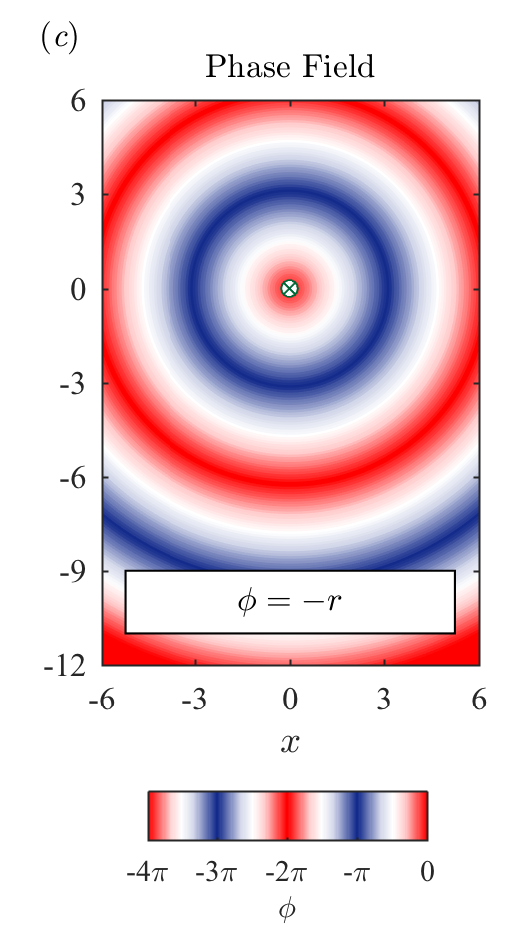}
\includegraphics{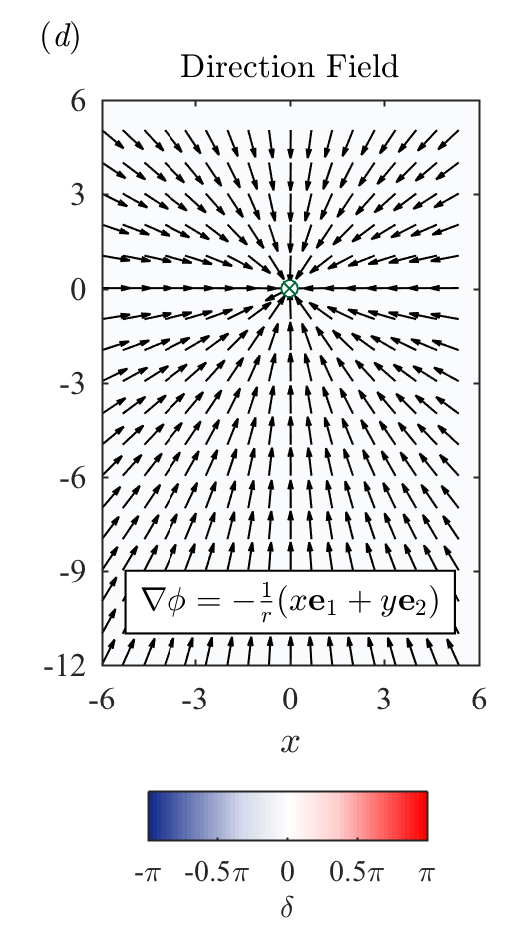}
\caption{%
Radially-symmetric signal field in~\eqref{eq:signal2} for $\ell = 6.5$.
(\textit{a})~Snapshot of the signal field $f(\mathbf{x},t)$ at time $t=0$,
(\textit{b})~the magnitude field $m(\mathbf{x})$,
(\textit{c})~the phase field $\phi(\mathbf{x})$,
(\textit{d})~the direction field $\nabla \phi$. Note that the associated alignment error $\delta(\mathbf{x})$ is identically zero.}
\label{fig:signal}
\end{figure*}

We consider a simplified signal field that locally approximates the signal generated by an oscillating source.
In particular, we consider a scalar signal field $f(r,\eta,t)\equiv f(r,t)$ propagating in a radially symmetric fashion (independent of $\eta$) from the origin,
\begin{equation}
 f(r,t) = 2\alpha \exp \left(\dfrac{-r}{\ell} \right) \cos \left( k r - \omega t \right).
\label{eq:signal}
\end{equation}
Here, $\alpha$, $k$ and $\omega$ are the signal magnitude, wavenumber and  frequency, respectively, and $\ell$ is a length scale that determines the spatial decay of the signal strength.

To simplify our analysis we use the time scale $1/\omega$, the length scale $1/k$, and the signal intensity scale $\alpha$ to define the dimensionless quantities
\begin{align*}
x^\dag = xk, \hspace{5pt}
y^\dag = yk, \hspace{5pt}
r^\dag = rk,  \hspace{5pt}
t^\dag = t\omega, \hspace{5pt}
\ell^\dag = \ell k, \hspace{5pt}
f^\dag = \dfrac{f}{\alpha}.
\end{align*}
For clarity, we drop the $(\cdot)^\dag$ notation and rewrite~\eqref{eq:signal} in nondimensional form,
\begin{equation}
f(r,t) = 2\exp \left( \dfrac{-r}{\ell} \right) \cos (r - t).
\label{eq:signal2}
\end{equation}
A snapshot of this signal is shown in figure~\ref{fig:signal}(\textit{a}).

Using~\eqref{eq:spectrum_definition},
the spectrum of~\eqref{eq:signal2} evaluated at $\omega = 1$ is given by
\begin{equation}
\hat{f}(r) =  \exp \left(-\dfrac{r}{\ell}-\i r \right).
\label{eq:spectrum}
\end{equation}
Substituting into~\eqref{eq:mag_phase_definition}, one gets the associated magnitude and phase fields,
\begin{equation}
m (r) = \exp \left(-\dfrac{r}{\ell}\right) \quad \quad \text{and} \quad \quad \phi(r) = -r.
\end{equation}
These fields are shown in figures~\ref{fig:signal}(\textit{b}) and (\textit{c}), respectively.
The normalized gradient of $\phi(r)$ is given by ${\nabla\phi}/{\|\nabla\phi\|} =  \mathbf{c}_1$ and the alignment error is identically zero, $\delta = 0$, everywhere in the $(r,\eta)$-plane, meaning that the spectral direction field 
$\nabla \phi$ always points towards the source; the normalized gradient field and contours of the alignment error are depicted in figure~\ref{fig:signal}(\textit{d}).

We substitute $\delta = 0$ into~\eqref{eq:eom_cframe} to get the equations of motion of the sensor in the simplified signal field $f(r,t)$,
\begin{equation}
\begin{split}
\dot{r} = &- V \cos \psi,\qquad
r\dot{\eta} = {V}\sin \psi,\\
r\dot{\psi} = & \left( V  - rG(r) \right)\sin\psi,
\end{split}
\label{eq:dynamics_radial_simplified}
\end{equation}
with the understanding that $V$ and $G$ are dimensionless.  
Due to the rotational symmetry of the signal field, the equations of motion are invariant under superimposed rotations on the field-sensor system. As a result, the equations for $r$ and $\psi$ in~\eqref{eq:dynamics_radial_simplified} decouple from the equation for $\eta$ and can be solved independently. The results can then be substituted into the equation $r\dot{\eta} = V\sin\psi$ to reconstruct $\eta(t)=\int_{0}^t {\rm d}\tau V\sin(\psi(\tau))/r(\tau)$. The two coupled equations in $(r,\psi)$ can be reduced further by noting that the phase gradient $\nabla \phi$ associated with the signal field in~\eqref{eq:signal2} is invariant under translations from $r$ to $r-t$,
which gives rise to an integral of motion $\mathcal{Q}$. To compute $\mathcal{Q}$, it is convenient to rewrite the equations for $r$ and $\psi$ in differential form,
\begin{equation}
\begin{split}
\dfrac{r}{V - rG} \dfrac{{\rm d}\psi}{\sin\psi} = -\dfrac{{\rm d}r}{V \cos\psi} ,
\end{split}
\end{equation}
which after straightforward manipulation yields the exact differential,
\begin{equation}
\begin{split}
{\rm d} \mathcal{Q} =  \dfrac{{\rm d}(\sin\psi)}{\sin\psi} + \left( \dfrac{1}{r} - \dfrac{G}{V}  \right){\rm d}r = 0.
\end{split}
\end{equation}
A direct integration of the above differential gives $\mathcal{Q} = r  \sin \psi \exp \left[ - \int {G(r)dr }/{V}  \, \right]=$ constant. Without loss of generality, we can scale $\mathcal{Q}$ by the constant $1/\rho = G_o/V$ to get
\begin{equation}
\mathcal{Q}(r,\psi) = \dfrac{r}{\rho} \sin \psi \exp \left[ - \int \dfrac{G(r)}{V}  \, dr \right] = {\rm constant}.
\label{eq:q_definition}
\end{equation}
We solve for $\psi$ in terms of $\mathcal{Q}$ and $r$ and substitute back into the equation for $r$ in~\eqref{eq:dynamics_radial_simplified} to get
\begin{equation}
\displaystyle
\dot{r} = \mp V \sqrt{1-\dfrac{\mathcal{Q}^2}{\dfrac{r^2}{\rho^2}  \exp \left( - \dfrac{2}{V}\int{G(r)} \, dr\right)}} \ .
\label{eq:radial_dyn}
\end{equation}
This equation can be readily integrated to get $r(t)$. 

\begin{figure*}
\centering
\includegraphics{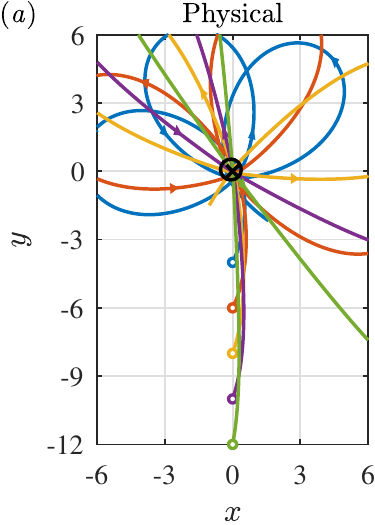}
\includegraphics{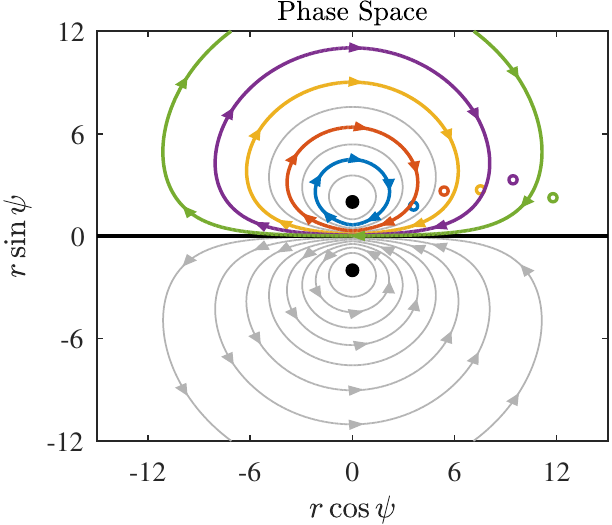}
\includegraphics{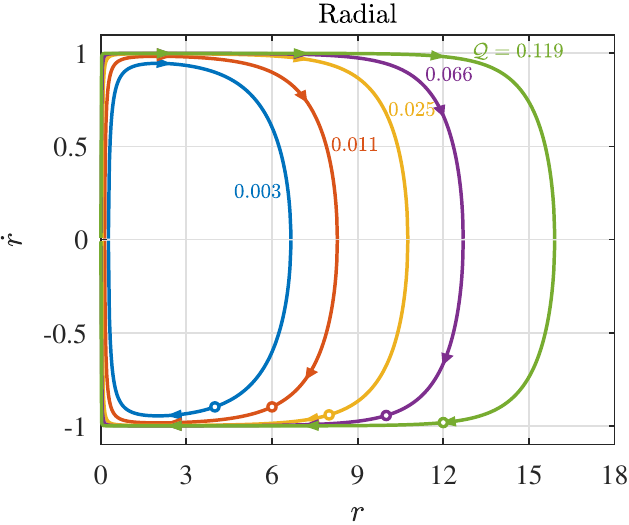}\\
\includegraphics{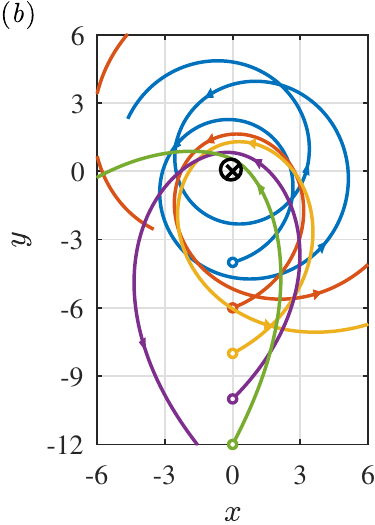}
\includegraphics{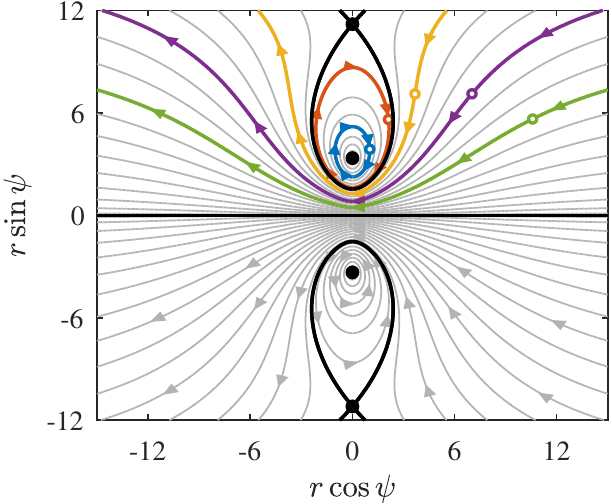}
\includegraphics{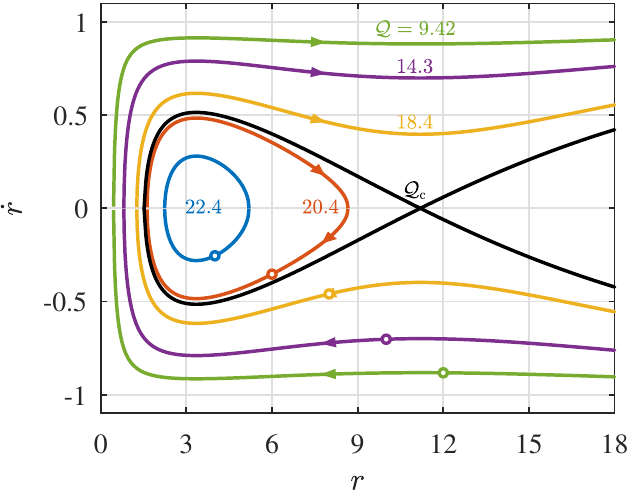}\\
\includegraphics{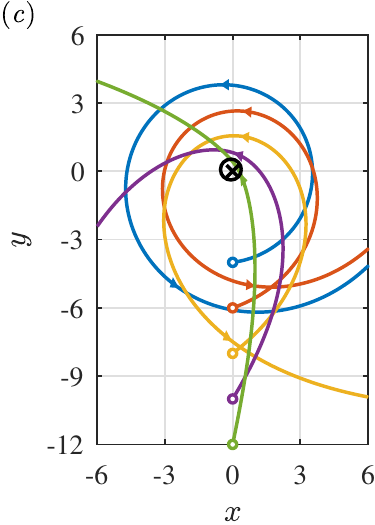}
\includegraphics{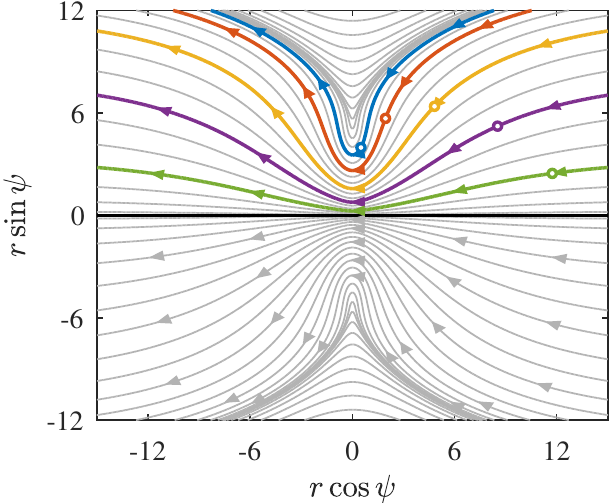}
\includegraphics{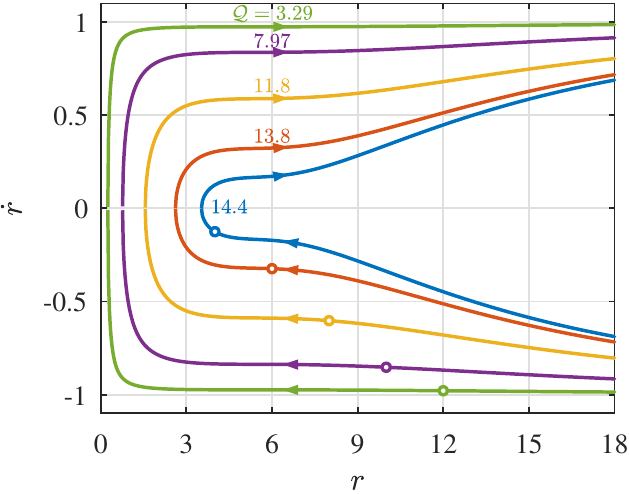}\\
\includegraphics{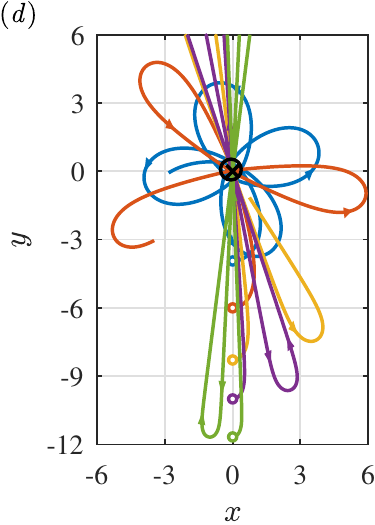}
\includegraphics{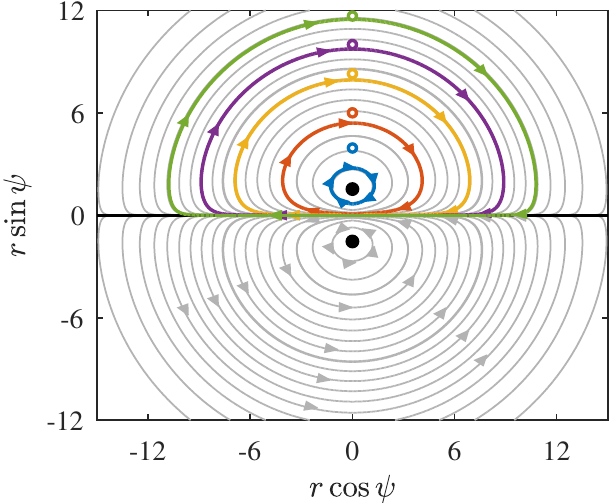}
\includegraphics{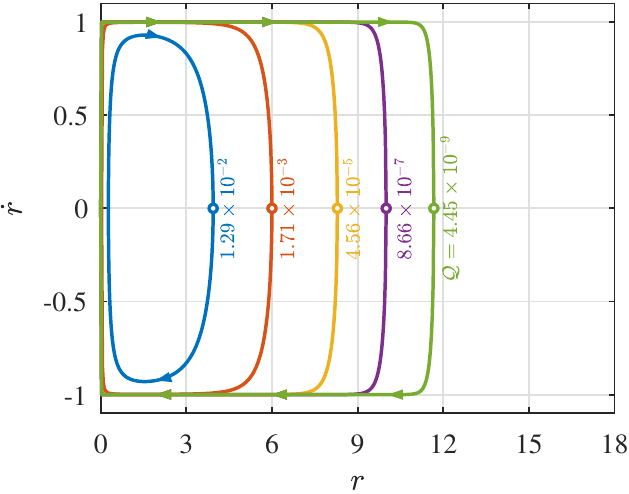}\\
\caption{%
Behavior of the mobile sensor in the radially-symmetric signal field in~\eqref{eq:signal2}. The left panel depicts trajectories in the $(x,y)$ plane starting at positions located successively further away from the source.
In the middle panel, trajectories are shown in the phase space $(r\cos\psi,r\sin\psi)$, with fixed points depicted as black dots.
Relative equilibria along $\sin \psi = 0$ are depicted as solid black lines.
The right panel depicts the solution curves in the $(r,\dot{r})$ space for various values of the integral of motion $\mathcal{Q}$.
For all cases $V = 1$ and $\rho = 2.0$.
(\textit{a})~Static gain. Two center-type fixed points occur. Unconditional convergence is achieved for all initial conditions.
(\textit{b})~Proportional gain with $\ell = 6.5$. Two center-type and two saddle-type fixed points occur.
The homoclinic orbits associated with the saddle fixed points are highlighted in black. Conditional convergence is achieved for initial conditions inside the closed-loop of these orbits.
(\textit{c})~Proportional gain with $\ell = 5.4$. No fixed points occur and no convergence.
(\textit{d})~Inverse gain with $\ell = 5.8$. Two center-type fixed points occur. Unconditional convergence is achieved for all initial conditions.
}
\label{fig:analytic_combined}
\end{figure*}

\section{Closed-Loop  Performance}
\label{sec:convergence}

We analyze the closed-loop performance of the mobile sensor in the simplified signal field in the context of three types of gain functions: (i) static gain of constant value $G_o$, (ii)  gain proportional to the  magnitude of the signal $G = G_o m(r)$, and (iii)  inversely-proportional gain  $G = G_o/m(r)$. In each case, we show the behavior of the sensor in the $(x,y)$ plane by numerically integrating the nonlinear equations of motion expressed in cartesian coordinates in~\eqref{eq:eom_scalar}. We then analyze the behavior in the phase space $(r\cos\psi,r\sin\psi)$ associated with the system of equations expressed in~\eqref{eq:dynamics_radial_simplified} after eliminating the $\eta$ equation. Lastly, we determine the success of the sensory-control law in following the signal field and \emph{converging} towards its source.

\subsection{Static Gain}
\label{sec:static_gain}

Consider the case of a static gain and set $G$ to a constant value $G_o$.
We solve the system of nonlinear equations in~\eqref{eq:eom_scalar} numerically for $V = 1$ and $G_o = 0.5$.
Figure \ref{fig:analytic_combined}\textit{a} (left) shows representative trajectories in the $(x,y)$ plane for distinct initial conditions, at progressively larger distances from the source, marked as `$\circ$' on the plot.
In each case, the mobile sensor approaches and loops around the target, providing multiple opportunities to intercept the source.

We analyze this behavior in the phase space $(r,\psi)$ where, from~\eqref{eq:dynamics_radial_simplified}, one has  
\begin{equation}
\begin{split}
\dot{r} &= -V\cos \psi, \quad \dot{\psi} = \left( \frac{1}{r} -\frac{1}{\rho} \right) V \sin \psi .
\end{split}
\label{eq:dynamics_static}
\end{equation}
This system of equation admits two \emph{fixed points} $(r^\ast, \psi^\ast)$ for which $\dot{r} = \dot{\psi} = 0$ for all time,
\begin{equation}
r^\ast = \rho, \quad \psi^\ast = \pm \frac{\pi}{2}.
\end{equation}
These fixed points are marked as `$\bullet$' in figure~\ref{fig:analytic_combined}a (middle) in the phase space 
$(r\cos\phi, r\sin\psi)$.
We determine the linear stability of these equilibria by computing the Jacobian matrix $J$ of~\eqref{eq:dynamics_static} and evaluating it at $(r^\ast, \psi^\ast)$
\begin{equation}
J\left(\rho, \pm \frac{\pi}{2} \right) = \pm V \begin{bmatrix}
0 & 1\\
-\rho^2 & 0
\end{bmatrix}.
\end{equation}
The corresponding eigenvalues $\lambda= \pm \i G_o$ are purely imaginary; thus, the fixed points are \emph{centers} by nature and are \emph{marginally stable}.
We also observe a line of degenerate \emph{relative} equilibria along $\sin \psi = 0$, for which $\psi = 0$ mod$(\pi)$ and $\dot{\psi} = 0$ and $\dot{r} = V$ for all time. Along this line, the sensor heads directly towards the source. 
We mark this line in figure~\ref{fig:analytic_combined}a (middle) and note that it separates the sensor trajectories that orbit clockwise around the source from those that orbit in a counter-clockwise fashion.

To complete the \emph{phase portrait} of the mobile agent, we superimpose onto figure 4a (middle) the  contours of the integral of motion $\mathcal{Q}$. For static gain, one has
\begin{equation}
\mathcal{Q} = \frac{r}{\rho} \sin \psi \exp\left( -\dfrac{r}{\rho} \right) = {\rm constant}.
\end{equation}
Generic contour lines are shown in gray. The specific contours corresponding to the trajectories shown in figure~\ref{fig:analytic_combined}a (left) are highlighted in color.

Lastly, in figure~\ref{fig:analytic_combined}\textit{a}~(right), we plot the radial phase velocity $\dot{r}$ as a function of $r$ for the five different values of $\mathcal{Q}$ corresponding to the trajectories in~\ref{fig:analytic_combined}\textit{a}~(left). Here, we use the radial dynamics in~\eqref{eq:radial_dyn}, which for static gain becomes
\begin{equation}
\label{eq:rdot_static}
\dot{r} = \pm V \sqrt{ 1-  \dfrac{\mathcal{Q}^2}{\dfrac{r^2}{\rho^2} \exp \left( -\dfrac{2r}{\rho} \right)} }.
\end{equation}
In each case we mark the initial condition and phase velocity by `$\circ$'.
We also label each curve by its corresponding value of $\mathcal{Q}$.
Clearly, for each value of $\mathcal{Q}$, the radial excursion is bounded from below and above.
To prove this for all initial conditions, we note that
physically-speaking, since $r$ is a positive scalar, equation~\eqref{eq:rdot_static} is well-defined only when the right-hand side is real-valued, that is to say, 
\begin{equation}
\left | \frac{\mathcal{Q}}{\dfrac{r}{\rho} \exp \left(-\dfrac{r}{\rho}\right)} \right | \leq 1.
\end{equation}
This condition can be used to determine the upper and lower bounds on the distance $r$,
\begin{equation}
- W_0 ( - |\mathcal{Q}|) \leq \frac{r}{\rho} \leq - W_{-1} ( -|\mathcal{Q}|),
\end{equation}
where $W_p(\cdot)$ is the $p^\text{th}$ Lambert product logarithm function.
Thus, we have shown that $r$ is bounded from both below and above for all values of $\mathcal{Q}$; in other words, $r$ remains finite for all time and all initial conditions.
The implication of this result is that regardless of the initial conditions, the trajectory of the mobile sensor  will circle around the source within a bounded distance $r$ from the source.
These periodic loops around the source have the practical advantage that the sensor will have multiple opportunities to approach the source.  We say that the mobile sensor exhibit  \emph{unconditional convergence} to the source.

\subsection{Proportional Gain}

We now consider a gain strategy that is proportional to the magnitude of the signal field $G = G_o \exp (-r/\ell)$.
Figure \ref{fig:analytic_combined}\textit{b} (left) shows the behavior of the mobile sensor in the $(x,y)$ plane for  $V = 1$, $G_o = 0.5$, and $\ell = 6.5$ and same initial conditions as in figure~\ref{fig:analytic_combined}(\textit{a}). 
Unlike in figure~\ref{fig:analytic_combined}\textit{a} (left), here, three trajectories (shown in green, purple, and yellow) get closer to the source once and then travel off to infinity, whereas two trajectories (shown in red and blue) make multiple loops around the source.
This implies a change in the stability characteristics of the system as the gain changes, with some initial conditions resulting in unstable trajectories, drifting off to $r\rightarrow \infty$ as time progresses.
To understand this change in behavior, we examine the dynamics in the phase space $(r\cos\psi,r\sin\psi)$, given by the equations of motion
\begin{equation}
\begin{split}
\dot{r} &= -V\cos \psi,\quad
\dot{\psi} = \left[ \frac{1}{r} -\dfrac{1}{\rho}\exp\left(-\dfrac{r}{\ell}\right) \right]V \sin \psi.
\end{split}
\label{eq:dynamics_proportional}
\end{equation}
Here, $\dot{\psi}$ is inversely proportional to an exponential function of $r$.
The physical interpretation is that, for this proportional gain, the mobile sensor will turn more sharply as it gets closer to the source. 

Equations~\eqref{eq:dynamics_proportional} admit four fixed points for which $\dot{r} = \dot{\psi} = 0$ for all time,
\begin{equation}
r^\ast = -\ell W_{-1} \left( - \frac{\rho}{\ell} \right), \quad \phi^\ast = \pm \frac{\pi}{2},
\end{equation}
and
\begin{equation}
r^{\ast\ast} = -\ell W_{0} \left( - \frac{\rho}{\ell} \right), \quad \phi^{\ast\ast} = \pm \frac{\pi}{2}.
\end{equation}
These fixed points are marked as `$\bullet$' in figure~\ref{fig:analytic_combined}b (middle).
To examine the linear stability of these fixed points, we evaluate the Jacobians
\begin{equation}
\begin{split}
J &= \pm V \begin{bmatrix}
0 & 1\\
\mathcal{X} & 0
\end{bmatrix},
\end{split}
\end{equation}
where 
\begin{equation}
\left. \mathcal{X} \right|_{\left( r^{\ast} , \phi^{\ast}\right)} = \dfrac{\mp 1}{\ell^2W_{-1}\left(-\dfrac{\rho}{\ell}\right)} \left( \dfrac{1}{W_{-1}\left(-\dfrac{\rho}{\ell}\right)}+1 \right),
\end{equation}
and
\begin{equation}
\left. \mathcal{X} \right|_{\left( r^{\ast\ast} , \phi^{\ast\ast}\right)} = \dfrac{\mp 1}{\ell^2W_{0}\left(-\dfrac{\rho}{\ell}\right)} \left( \dfrac{1}{W_{0}\left(-\dfrac{\rho}{\ell}\right)}+1 \right).
\end{equation}
This leads to the eigenvalue pairs
\begin{equation}
\begin{split}
\frac{\lambda^\ast \ell}{\sqrt{V}} &= \pm \sqrt{ \dfrac{1}{W_{-1}(-\dfrac{\rho}{\ell})}\left( \dfrac{1}{W_{-1}(-\dfrac{\rho}{\ell})} + 1 \right)},
\end{split}
\end{equation}
and
\begin{equation}
\begin{split}
\frac{\lambda^{\ast\ast} \ell}{\sqrt{V}} &= \pm \i \sqrt{ \dfrac{1}{W_{0}(-\dfrac{\rho}{\ell})}\left( \dfrac{1}{W_{0}(-\dfrac{\rho}{\ell})} + 1 \right)},
\end{split}
\end{equation}
for $\left( r^{\ast} , \phi^{\ast}\right)$ and $\left( r^{\ast\ast} , \phi^{\ast\ast}\right)$, respectively.
The fixed points $(r^\ast, \phi^\ast)$ are unstable (saddle type) while $(r^{\ast \ast}, \phi^{\ast \ast})$ are linearly stable (centers).
As before, $ \psi = 0$ mod $(\pi)$ corresponds to a degenerate family of relative equilibria.

We complete the phase portrait in figure~\ref{fig:analytic_combined}\textit{b} (middle) by plotting the contour lines of $\mathcal{Q}$ given by 
\begin{equation}
\mathcal{Q} = \dfrac{r}{\rho} \sin \psi \exp \left[\dfrac{\ell}{\rho} \exp\left( - \dfrac{r}{\rho}\right)\right].
\end{equation}
Notice the  \emph{homoclinic-type separatrices}, highlighted with a thickened black curve, associated with the saddle fixed points.
Trajectories that start inside the homoclinic orbit are stable and those outside are unstable. 
This distinction explains the difference in behavior we observed in figure~\ref{fig:analytic_combined}\textit{b} (left).

In figure~\ref{fig:analytic_combined}\textit{b}~(right), we plot the radial phase velocity $\dot{r}$ as a function of $r$ for the five different values of $\mathcal{Q}$ corresponding to the trajectories in~\ref{fig:analytic_combined}\textit{b}~(left).
Here, we use the radial dynamics in~\eqref{eq:radial_dyn}, which for proportional gain becomes
\begin{equation}
\dot{r} = \pm V \sqrt{1-  \frac{\mathcal{Q}^2}{\dfrac{r^2}{\rho^2} \exp\left(2\dfrac{\ell}{\rho} \exp\left(-\dfrac{r}{\ell}\right)\right)}}\ .
\label{eq:proportional_radial}
\end{equation}
We define the \emph{critical conserved quantity} $\mathcal{Q}_\text{cr}$ as the value of $\mathcal{Q}$ corresponding to the separatrices associated with the saddle points $(r^{\ast},\pm \pi/2)$. That is to say, we evaluate $\mathcal{Q}$ at $(r^{\ast},\pm \pi/2)$ to get
\begin{equation}
\mathcal{Q}_\text{cr} = \mp \dfrac{\ell}{\rho} W_{-1} \left( - \frac{\rho}{\ell} \right) \exp \left[ \frac{-1}{W_{-1} \left(-\dfrac{\rho}{\ell}\right)}\right].
\end{equation}
If $| \mathcal{Q} | < | \mathcal{Q}_{\text{cr}} |$, the sensor will approach the source only once, as seen in the green, purple, and yellow curves.
On the other hand, if $| \mathcal{Q} | > | \mathcal{Q}_{\text{cr}} |$, the trajectory will either  loop around the source indefinitely or will diverge to infinity depending on whether the initial radial position is greater than or less than $r^\ast$, respectively. 
Given that the radial excursion is bounded from above for only a subset of the initial conditions, we will call this \emph{conditional convergence}.

Let us consider the case where $\ell = 5.4$.
Keeping $V$ and $G_o$ the same, we solve the system of nonlinear equations in~\eqref{eq:kinematics} numerically, generating the trajectories in the $(x,y)$ plane which we plot in figure \ref{fig:analytic_combined}\textit{c} (left).
Clearly, all of the the trajectories make one close-in pass of the source and then move away to infinity.

To understand this behavior, we reexamine~\eqref{eq:dynamics_proportional} with the lower value of $\ell$ and find that there are no fixed points.
This is because a \emph{saddle-node bifuraction} occurs when $\ell = \rho e$ and the saddle- and center-type fixed points collide and \emph{annihilate}.
Upon further examination of figure \ref{fig:analytic_combined}\textit{c} (middle), we confirm that all trajectories regardless of initial condition are unbounded.
Finally, the radial dynamics in figure~\ref{fig:analytic_combined}\textit{c} (right) confirm that $r$ is unbounded for all initial conditions. The sensor is said to \emph{diverge}.

\subsection{Inversely Proportional Gain}

Lastly, we consider the gain strategy $G = G_o \exp (r/\ell)$ that is inversely proportional to the magnitude field. That is, the gain decreases with signal strength.
In figure~\ref{fig:analytic_combined}\textit{d} (left), we plot the motion of the sensor in the $(x,y)$ plane for  $V = 1$, $G_o = 0.5$, and $\ell = 5.8$. In this case, the trajectories resemble those in figure~\ref{fig:analytic_combined}\textit{a} (left), but exhibit longer loops with tighter turns around the source.
As before, we examine the dynamics expressed in the phase space $(r\cos\psi,r\sin\psi)$, 
\begin{equation}
\begin{split}
\dot{r} &= -V \cos \psi, \quad \dot{\psi} =  \left[ \frac{1}{r} -\dfrac{1}{\rho}\exp\left(\dfrac{r}{\ell}\right) \right] V \sin \psi.
\end{split}
\label{eq:dynamics_inverse}
\end{equation}
Here, the rate of change of the heading angle $\dot{\psi}$ is proportional to  $\exp(r/\ell)$, implying that the mobile sensor turns more sharply when it is further away from the source. This is in contrast to the proportional gain strategy $\dot{\psi} \sim \exp(-r/\ell)$, where the gain increases with signal strength and the sensor turns more sharply closer to the source.

Equations~\eqref{eq:dynamics_inverse} admit two fixed points for which $\dot{r} = \dot{\psi} = 0$ for all time,
\begin{equation}
r^\ast = \ell W_0\left( \frac{\rho}{\ell} \right), \quad \quad \phi^\ast = \pm \frac{\pi}{2},
\end{equation}
which are depicted in figure~\ref{fig:analytic_combined}\textit{d} (middle) as `$\bullet$'.
To examine the linear stability of these fixed points, we evaluate the Jacobian
\begin{equation}
J \left(r^\ast , \pm \frac{\pi}{2} \right) = \pm V \begin{bmatrix}
0 & 1\\
\mathcal{X} & 0
\end{bmatrix},
\end{equation}
where
\begin{equation}
\mathcal{X} = \frac{\mp 1}{\ell^2 W_0\left( \dfrac{\rho}{\ell}\right)} \left( \frac{1}{W_0\left(\dfrac{\rho}{\ell}\right)} + 1 \right).
\end{equation}
This leads to the eigenvalue pairs
\begin{equation}
\frac{\lambda \ell}{\sqrt{V}} = \pm \i \sqrt{ \frac{1}{W_0\left(\dfrac{\rho}{\ell}\right)} \left( \frac{1}{W_0\left(\dfrac{\rho}{\ell}\right)} + 1\right)}.
\end{equation}
Since these eigenvalue pairs are purely imaginary and equal in magnitude,  the fixed point is a center and is linearly stable.

We complete the phase portrait in  figure~\ref{fig:analytic_combined}\textit{d} (middle) by plotting the contour lines of the conserved quantity,
\begin{equation}
\mathcal{Q} = \frac{r}{\rho} \sin \psi \exp \left[ -\frac{\ell}{\rho} \exp \left( \frac{r}{\ell} \right) \right],
\end{equation}
The contour lines are closed orbits around the fixed points, in agreement with the fact that these points are of center type.
Importantly, all orbits are periodic and bounded, meaning that the mobile sensor does not travel off towards $r\rightarrow \infty$.

The boundedness of the radial distance from the source $r$ can be best seen by examining the radial equation of motion
\begin{equation}
\dot{r} = \pm V \sqrt{ 1-  \frac{\mathcal{Q}^2}{\dfrac{r^2}{\rho^2} \exp\left(\dfrac{2r}{\ell}\right)} },
\label{eq:inverse_radial}
\end{equation}
We plot the radial velocity  $\dot{r}$ versus $r$ in figure~\ref{fig:analytic_combined}\textit{d} (right) for all five different trajectories highlighted in figure~\ref{fig:analytic_combined}\textit{d} (left).
The radial excursion $r$ is bounded from both below and above for all initial conditions and hence the sensor dynamics is said to be  \emph{unconditionally convergent}.

\section{Locating an Oscillating Airfoil by Following its Wake}

\begin{figure*}
\centering
\includegraphics{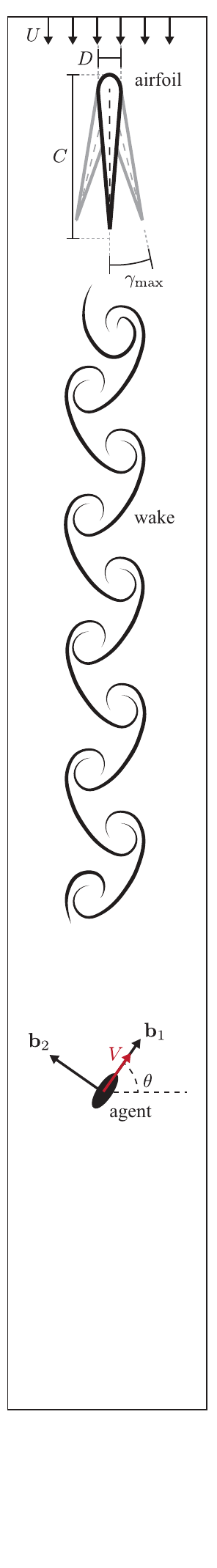}
\includegraphics{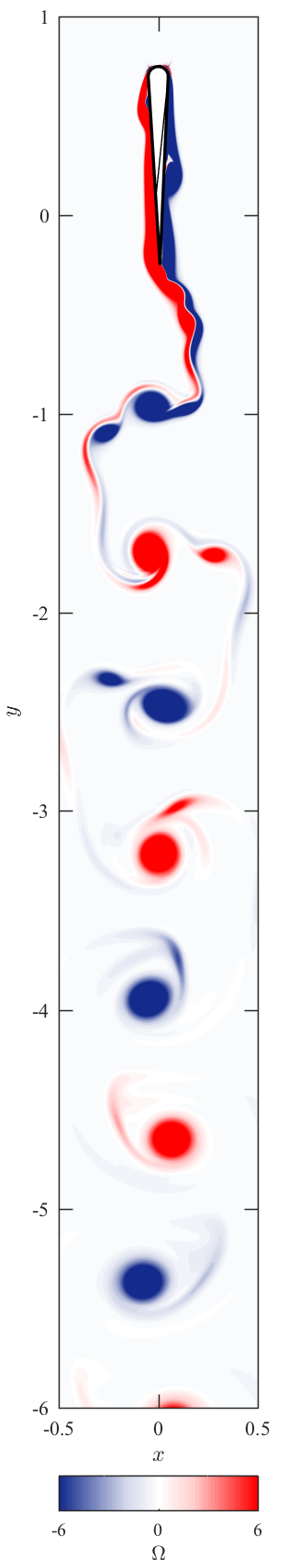}
\includegraphics{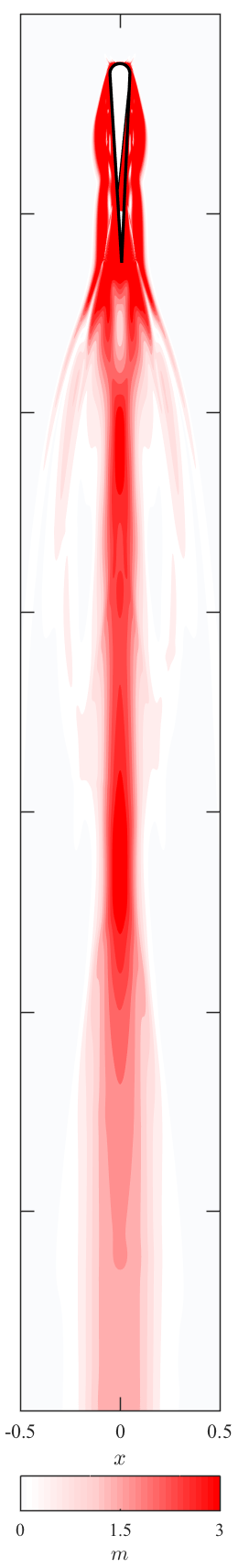}
\includegraphics{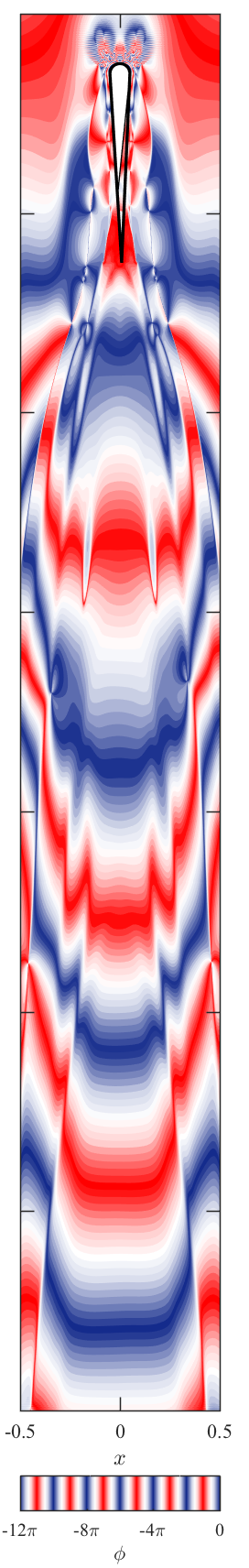}
\includegraphics{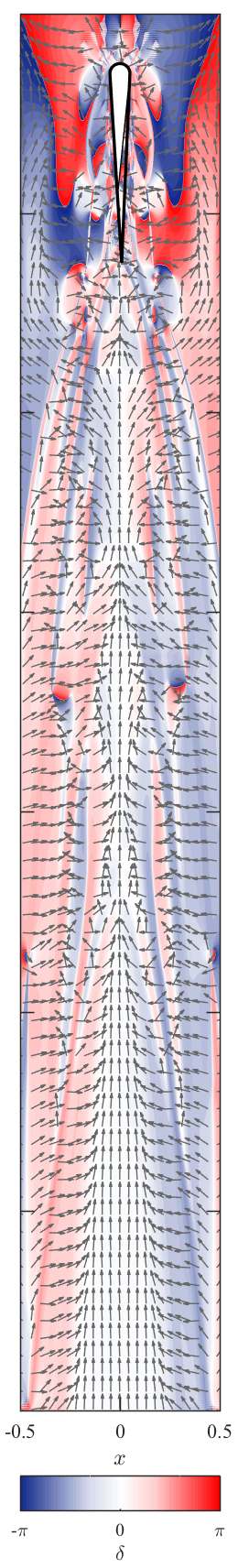}
\caption{%
Behavior of the mobile sensor in the wake of an oscillating airfoil in a uniform incoming flow of strength $U$. 
The airfoil oscillates at dimensionless frequency St (Strouhal number) and amplitude A.
A wake forms downstream of the oscillating airfoil, with regions of spatiotemporally-correlated coherent vortical structures.
(\textit{a})~Snapshot of vorticity field field $\Omega(\mathbf{x},t)$ for $\text{Re} = 500$, $\text{St} = 0.075$ and $\text{A} = 1.66$.
(\textit{b})~the magnitude field $m(\mathbf{x})$ of the spectral vorticity,
(\textit{c})~the phase field $\phi(\mathbf{x})$,
(\textit{d})~the direction field $\nabla \phi / \| \nabla \phi \|$ and the associated alignment error $\delta(\mathbf{x})$.
}
\label{fig:wake_fields}
\end{figure*}

Now that we analyzed the convergence properties of the mobile sensor in a simplified signal field, we seek to demonstrate its performance in a more complex signal field arising from the flow past an oscillating airfoil.
To this end, we consider an airfoil-shaped body with diameter $D$ and chord $C$ in an otherwise uniform flow of speed $U$.
We control the pitch angle $\gamma(t)$ sinusoidally with frequency $f$ such that
$\gamma(t) = \gamma_\text{max} \cos(2\pi f t)$,
as depicted in figure~\ref{fig:wake_fields}(\textit{a}).
The parameters can be reduced to three dimensionless parameters: the {Strouhal number}  $\text{St} = fD/U$, the {amplitude ratio} $\text{A} = 2C \sin \gamma_\text{max}/D$, and the Reynolds number $\text{Re} =  U D/\nu$, where $\nu$ is the fluid viscosity. 
These dimensionless parameters characterize the structure of the resultant wake; see, e.g.,~\cite{Schnipper2009,Colvert2018}.

The flow around the airfoil is described by the fluid velocity vector field $\mathbf{u}(\mathbf{x},t)$ and the scalar pressure field $p(\mathbf{x},t)$.
The time evolution of these fields is governed by the incompressible Navier-Stokes equations
\begin{equation}
\begin{split}
&\frac{\partial \mathbf{u}}{\partial t} = -\mathbf{u} \cdot \nabla \mathbf{u} -\nabla p + \frac{1}{\text{Re}} \nabla^2 \mathbf{u},\\
&\nabla \cdot \mathbf{u} = 0.
\end{split}
\end{equation}
together with the no-slip condition at the airfoil boundary. 
We solve these equations numerically for $\text{Re} = 500$ and $\text{St} = 0.075$ using an 
algorithm based on the Immersed Boundary Method~\cite{Peskin1977}. The algorithm is described in detail in~\cite{Mittal2008}, and has been implemented, optimized and tested extensively by the group of Haibo Dong; see, e.g.,~\cite{Vargas2008, Bozkurttas2009}. 

We calculate the vorticity field $w(\mathbf{x},t) = (\nabla \times \mathbf{u} ) \cdot \mathbf{e}_3$ where $\mathbf{e}_3 = \mathbf{e}_1 \times \mathbf{e}_2$ is the unit vector normal to the plane.
Figure~\ref{fig:wake_fields}(\textit{a}) shows contours of vorticity for a snapshot in time.
We take $w(\mathbf{x},t)$ as the signal field.
 We compute its Fourier transform $\hat{w}(\mathbf{x},\omega_n)$, evaluated at the first Fourier mode $n=1$, and calculate its magnitude $m(\mathbf{x})$ and phase $\phi(\mathbf{x})$, as depicted in figures~\ref{fig:wake_fields}(\textit{b}) and (\textit{c}), respectively.
The direction field $\nabla \phi/ \| \nabla \phi\|$ and the alignment error $\delta(\mathbf{x})$ are plotted in figure~\ref{fig:wake_fields}(\textit{d}).

Although the structure of the vorticity signal field $w(\mathbf{x},t)$ and its spectral fields in figure~\ref{fig:wake_fields} is far more complex than the radially-symmetric signal given in~\eqref{eq:signal2} and plotted in figure~\ref{fig:signal}, the two signals have some similarities.
First, the magnitude $m(\mathbf{x})$ of the signal increases as we approach the airfoil.
Second, the direction field $\nabla \phi/ \| \nabla \phi\|$ points in the upstream direction over a large area downstream of the airfoil.
These two components lead us to believe that the control law proposed in this work will drive a mobile sensor initially positioned downstream to follow the hydrodynamic signal and converge to the airfoil. 
In fact, downstream of the airfoil, the alignment error $\delta(\mathbf{x})$ is nonzero; to the left of the centerline, the alignment error is largely positive, and the  vector field $\nabla \phi/ \| \nabla \phi\|$ points towards the centerline.
To the left of the centerline, the alignment error is largely negative, and $\nabla \phi/ \| \nabla \phi\|$ points towards the centerline.
Therefore, a mobile sensor initially positioned in the wake of the airfoil and re-orienting according to $\nabla \phi/ \| \nabla \phi\|$ will be directed towards the source.
Following the analysis in section~\ref{sec:convergence}, we test this prediction in three types of gains: static, proportional and inversely-proportional. 

\begin{figure*}
\centering
\includegraphics{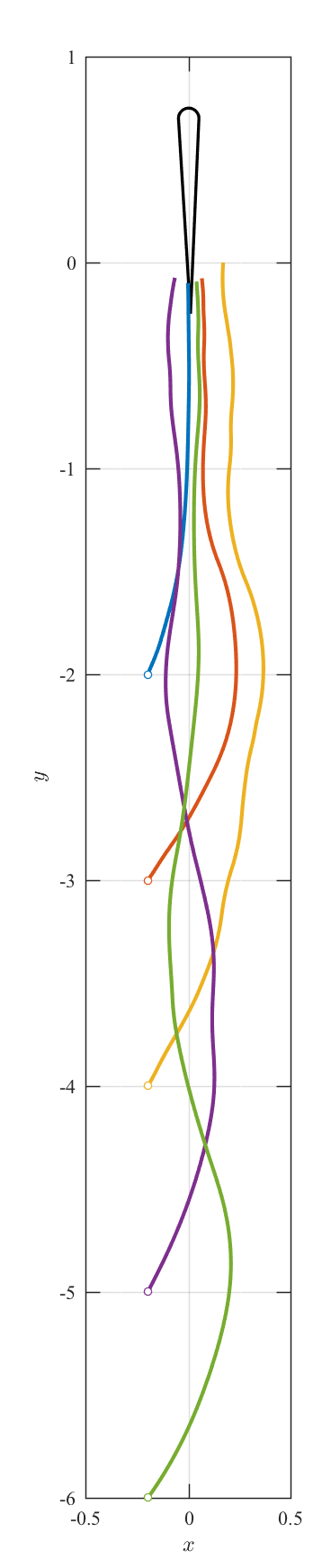}
\includegraphics{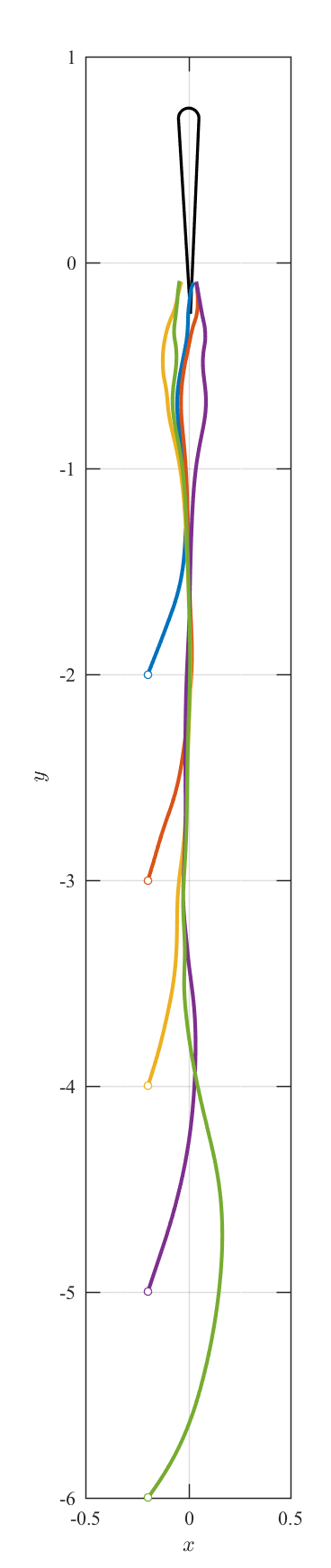}
\includegraphics{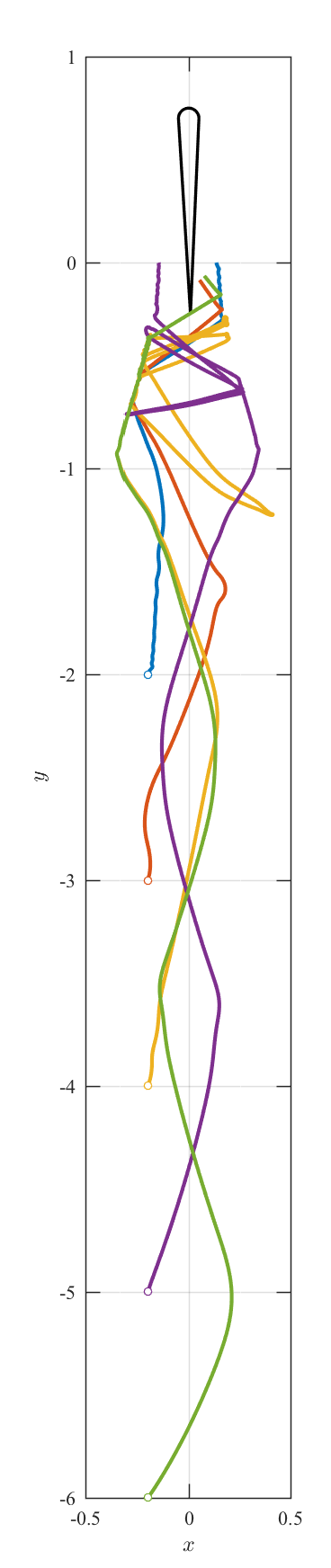}
\caption{%
Trajectories of the mobile sensor in the wake of the oscillating airfoil for the parameter values shown in figure~\ref{fig:wake_fields}:
(\textit{a})~static gain law,
(\textit{b})~proportional gain law, and
(\textit{c})~inversely proportional gain law.
In all cases, the  parameters of the mobile sensor are set to $G_o = 1$ and $V=1$.}
\label{fig:wake_trajectories}
\end{figure*}

\paragraph{Static Gain}
We solve~\eqref{eq:eom_scalar} numerically using the signal from the vorticity field for various initial conditions located downstream of the airfoil at successively larger distances from the source.
Results are shown in figure~\ref{fig:wake_trajectories}(\textit{a}) in the $(x,y)$ plane.
The mobile sensor largely oscillates about the centerline until it reaches the airfoil.
For the trajectories shown in green, red, and blue, the mobile sensor ends up very close to the airfoil.
However, for the purple and yellow trajectories, the mobile sensor's final position is further off the centerline.
In particular, the yellow trajectory is an example of an overshoot, where the mobile sensor did not turn sharply enough to stay in the wake.
Overall, the trajectories  end up in the vicinity of the airfoil but do not concentrate along the centerline -- instead, they are distributed transversally to the wake.
This is due to the fact that the gain strategy does not take into consideration any information about the signal strength, which is largely concentrated along the centerline.

\paragraph{Proportional Gain}
We incorporate the signal strength by choosing a gain strategy  $G= G_o m(\mathbf{x})$ proportional to the signal magnitude.
We again solve~\eqref{eq:eom_scalar} for the same initial conditions. The resulting trajectories in the $(x,y)$ plane are plotted in
figure~\ref{fig:wake_trajectories}(\textit{b}).
For this choice of initial conditions, the mobile sensor approaches then tracks the centerline of the wake.
The proportional gain strategy causes the mobile sensor to turn more sharply in response to higher signal magnitude.
Since the magnitude of the signal is concentrated along the centerline of the wake, the mobile sensor stays more closely aligned to it once it reaches the centerline.
Therefore,  the trajectories track the airfoil and approach it very closely.

\paragraph{Inversely-Proportional Gain}
Lastly, we consider the gain strategy where $G$ is inversely proportional to the signal magnitude such that $G= G_o/m(\mathbf{x})$.
Figure~\ref{fig:wake_trajectories}(\textit{a}) shows the trajectories of the sensor for the same initial conditions considered in
figures~\ref{fig:wake_trajectories}(\textit{a,b}).
Clearly, the mobile sensor oscillates about the centerline until it reaches the airfoil. The sensor tends to undergo more lateral excursions as it approaches the airfoil.

In the regions of the wake where the magnitude field is vanishingly small, the gain function becomes singular and the control law is ill-posed. Therefore, this algorithm is not practical for online implementation given that in real applications, there is no guarantee that the signal field will be always non-zero.

\section{Conclusions}

We presented a novel framework for analyzing the information embedded in time-periodic flow fields in the context of a source-seeking problem, where a mobile sensor, with only local information, is required to track the flow and locate its generating source.
Our approach takes advantage of the ``traveling wave"-like character of the signal field in order to reorient the mobile sensor in the direction opposite to the direction of propagation of the signal. 
Specifically, we measure the direction of signal propagation in the frequency domain and map back to the time domain to reorient the sensor.
This method is best suited for systems where the time scale associated with the signal is much smaller than the time scale associated with the sensor. 
We rigorously analyzed the convergence of the sensor to the source in the context of a simplified signal field. 
Then, through carefully laid out simulations of fluid flow past an oscillating airfoil, we demonstrated the efficacy of this sensory-control system in tracking the flow signal and locating the airfoil.

A few comments on the limitations of the algorithm are in order.
The transformation of time domain information at one position in space to the frequency domain requires in principle that the mobile sensor remains stationary for one period of the signal. 
However, waiting to collect the signal may not be feasible or practical to enforce in certain applications. 
In this case, we would need to extend the algorithm to perform true real-time spectral estimation, for which there are a variety of useful algorithms in the spectral analysis literature (see e.g.~\cite{Allen1977, Kuo1993}).
Another limitation of the model is that the mobile sensor is one-way coupled to the signal field --- that is, it can probe the signal but it does not affect it.
This condition was used to simplify the analysis and demonstrate the effectiveness of the sensing scheme.
Further research is required to understand the coupling of flow sensing to self-generated flow disturbances as done in~\cite{Gao2013,Free2018,Yen2018}.

Our results, that the time history of local flow quantities contain information that can be used to track the flow and locate its source, can be seen as a first step towards laying a foundation for source seeking in fluid flows.
Indeed, in~\cite{Colvert2018}, we demonstrated, using neural networks, that local time histories of the flow contain information relevant to \emph{categorizing} the flow field and hence its source, whereas here we showed, using a deterministic control law, that local time histories contain information relevant to \emph{finding} the source.
Future extensions of these works will unify these two strategies, by combining tools from artificial intelligence and machine learning with rigorous dynamical systems and control theoretic methods,  in a robust manner to allow mobile ``agents" equipped with hydrodynamic sensing capabilities to \emph{identify} and \emph{localize} other objects or agents in the flow field.

\bibliographystyle{IEEEtran}
\bibliography{source_seeking.bib}

\bigskip\bigskip
\bigskip\bigskip
\bigskip\bigskip
\bigskip\bigskip
 \begin{wrapfigure}{l}{0.1\textwidth}
  \begin{center}
    \includegraphics[width=0.1\textwidth]{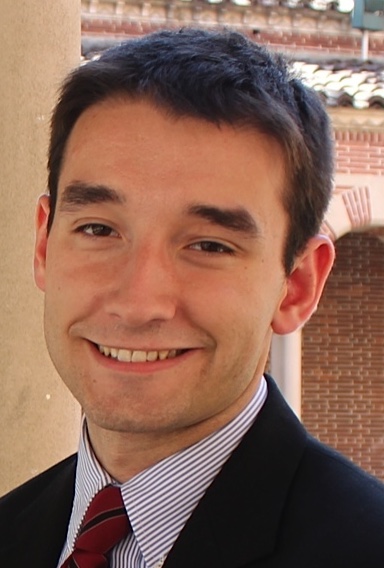}
     \end{center}
\end{wrapfigure}
\noindent \textbf{Brendan Colvert} received his B.S. ('14) and M.S. ('15) degrees in Aerospace Engineering from the University of Southern California (USC). He is currently working towards the Ph.D. degree under the direction of Dr. Eva Kanso at USC. The focus of his research is physics-based and data-driven modeling techniques for bio-inspired flow sensing and motion planning. Brendan Colvert is a National Defense Science and Engineering Graduate (NDSEG) Fellow.

\bigskip

 \begin{wrapfigure}{l}{0.1\textwidth}
  \begin{center}
    \includegraphics[width=0.1\textwidth]{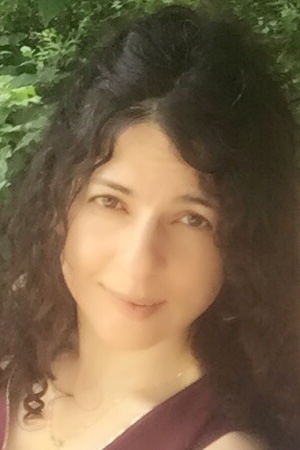}
     \end{center}
\end{wrapfigure}
\noindent
\textbf{Eva Kanso} is a professor and Z.H. Kaprielian Fellow in Aerospace and Mechanical Engineering at the University of Southern California. Prior to joining USC, Kanso held a two-year postdoctoral position (2003-2005) in Control and Dynamical Systems at Caltech. She received her Ph.D. (2003) and M.S. (1999) in Mechanical Engineering as well as her M.A. (2002) in Mathematics from UC Berkeley. Kanso is the recipient of an NSF early CAREER development award (2007), a Junior Distinguished Alumnus award from the American University of Beirut (2014), a USC Graduate Student Mentoring Award (2016), and an NSF INSPIRE grant (2017). Kanso's research focuses on the biophysics of living systems, with applications ranging from underwater and aerial locomotion, sensing and control to flow-mediated collective phenomena, biomedical and cilia-driven flows.

\end{document}